\documentclass[leqno,12pt]{article}

\usepackage{amssymb}
\usepackage{euscript}
\usepackage{flafter}
\usepackage{pstricks}
\usepackage{graphicx}
\usepackage{epstopdf} 

 \usepackage[latin1]{inputenc}

\usepackage{verbatim}

\usepackage{amsmath}

\usepackage{mathrsfs} 

\evensidemargin\oddsidemargin

\renewcommand{\thefootnote}{\fnsymbol{footnote}}

\begin{document}
\baselineskip=18pt
\setcounter{page}{1}

\renewcommand{\theequation}{\thesection.\arabic{equation}}
\newtheorem{theorem}{Theorem}[section]
\newtheorem{lemma}[theorem]{Lemma}
\newtheorem{definition}[theorem]{Definition}
\newtheorem{proposition}[theorem]{Proposition}
\newtheorem{corollary}[theorem]{Corollary}
\newtheorem{remark}[theorem]{Remark}
\newtheorem{fact}[theorem]{Fact}
\newtheorem{problem}[theorem]{Problem}
\newtheorem{conjecture}[theorem]{Conjecture}
\newtheorem{claim}[theorem]{Claim}

\newcommand{\eqnsection}{
\renewcommand{\theequation}{\thesection.\arabic{equation}}
    \makeatletter
    \csname  @addtoreset\endcsname{equation}{section}
    \makeatother}
\eqnsection

\def\r{{\mathbb R}}
\def\e{{\mathbb E}}
\def\p{{\mathbb P}}
\def\q{{\mathbb Q}}
\def\P{{\bf P}}
\def\E{{\bf E}}
\def\Q{{\bf Q}}
\def\z{{\mathbb Z}}
\def\N{{\mathbb N}}
\def\T{{\mathbb T}}
\def\G{{\mathbb G}}
\def\F{{\mathscr F}}
\def\L{{ \mathscr L}}
\def\A{{\mathscr A}}

\def\ee{\mathrm{e}}
\def\d{\, \mathrm{d}}
\def\law{{\buildrel \mbox{\small\rm (law)} \over =}}
\def\wcv{{\buildrel \mbox{\small\rm (law)} \over \longrightarrow}}


\vglue50pt

\centerline{\Large\bf  Local times of  subdiffusive  biased  walks  on trees} 

{
\let\thefootnote\relax\footnotetext{\scriptsize Partly supported by ANR project MEMEMO2 (2010-BLAN-0125).}
}

\bigskip
\bigskip


\medskip

 \centerline{Yueyun Hu\let\thefootnote\relax\footnote{\scriptsize LAGA, Universit\'e Paris XIII, 99 avenue J-B Cl\'ement, F-93430 Villetaneuse, France, {\tt yueyun@math.univ-paris13.fr}}}

\medskip

 \centerline{\it Universit\'e Paris XIII}

\bigskip

\bigskip
\bigskip

{\leftskip=2truecm \rightskip=2truecm \baselineskip=15pt \small

\noindent{\slshape\bfseries Summary.}  Consider a class of null-recurrent randomly biased walks on a super-critical Gaton-Watson tree.  We obtain the rates of convergence  of   the local times and  the quenched local probability for the biased walk in the sub-diffusive case.  These results are a consequence of  a sharp estimate on the return time, whose analysis is driven by a family of concave recursive equations on trees.

\bigskip


\bigskip


} 

\bigskip
\bigskip

\section{Introduction}
   \label{s:intro}


 We are interested in  a randomly biased walk $(X_n)_{n\ge0}$ on a supercritical Galton--Watson tree $\T$, rooted at  $\varnothing$. For any vertex $x\in \T \backslash \{ \varnothing\}$, denote by ${\buildrel \leftarrow \over x}$ its parent.    Let $\omega := (\omega(x, \cdot), \, x\in \T)$ be a sequence of vectors such that  for each vertex $x\in \T$,   $\omega(x, \, y) \ge 0$ for all $y\in \T$ and   $\sum_{y\in \T} \omega(x, \, y) =1$. We assume that  $\omega(x, \, y)>0$ if and only if either ${\buildrel \leftarrow \over x}=y$ or ${\buildrel \leftarrow \over y}=x$.

  For the sake of presentation, we add a specific vertex  ${\buildrel \leftarrow \over \varnothing}$, considered as the parent of $\varnothing$.  Let us stress that  ${\buildrel \leftarrow \over\varnothing} \not\in  \T$. We define $\omega({\buildrel \leftarrow \over\varnothing}, \varnothing):=1$ and modify the vector $\omega(\varnothing, \cdot)$ such that $\omega(\varnothing, {\buildrel \leftarrow \over \varnothing})>0$ and $\omega(\varnothing, {\buildrel \leftarrow \over \varnothing})+ \sum_{x: {\buildrel \leftarrow \over x}=\varnothing} \omega(\varnothing, x)=1$.

For given $\omega$,   $(X_n, \, n\ge 0)$ is a Markov chain  on $\T\cup\{ {\buildrel \leftarrow \over\varnothing}\}$ with transition probabilities $\omega$,   starting  from  $ \varnothing$;  i.e. $X_0= \varnothing$ and 
$$
P_\omega \big( X_{n+1} = y \, | \, X_n =x \big) 
= 
\omega(x, \, y).
$$


For any vertex $x\in \T $,  let $(x^{(1)}, \cdots, x^{(\nu_x)})$ be its children, where $\nu_x \ge 0$ is the number of children of $x$. Define ${\bf A}(x) := (A(x^{(i)}), \, 1\le i\le \nu_x)$ by
\begin{equation*}
    A(x^{(i)}) 
    := 
    \frac{\omega(x, \, x^{(i)})}{\omega(x, \, {\buildrel\leftarrow \over x})},
    \qquad 1\le i\le \nu_x \, .
\end{equation*}

When all $A(x^{(i)})=\lambda$ with some positive constant $\lambda$, the walk is called  $\lambda$-biased walk on a Galton-Watson tree and  was studied in detail by  Lyons, Pemantle and Peres   \cite{lyons-pemantle-peres-ergodic, lyons-pemantle-peres96}. We mention that   several conjectures  in \cite{lyons-pemantle-peres96} still remain open and we refer to Aidekon \cite{elie-vitesse}  for an explicit  formula on the speed of the  $\lambda$-biased walk and the references therein for recent developments.

When $A(x^{(i)})$  is also a random variable, the couple $(\T, \omega)$ is a marked tree in the sense of  Neveu \cite{neveu}, and the biased walk $X$  can be reviewed as a random walk in random environment.

Let us assume that ${\bf A}(x), x\in \T$ (including $x=\varnothing$) are i.i.d., and  denote  the vector  ${\bf A}(\varnothing)$ by $(A_1, ..., A_{\nu})$ for notational convenience. As such, $\nu\equiv \nu_\varnothing$ is the number of children of $\varnothing$.     Denote by  $\P$  the law of $\omega$ and define $\p(\cdot):=\int P_\omega(\cdot) \P(d\omega)$. In the literature  of random  walk in random environment, $P_\omega$ is referred to the quenched probability whereas $\p$ is  the annealed probability.

 Define $$  \psi(t):= \log \E \Big( \sum_{i=1}^\nu \, A_i^t\Big) \in (-\infty, \infty],   \qquad \forall t \in \r. $$
 
\noindent  In particular, $\psi(0)= \log \E(\nu) >0$ since $\T$ is supercritical.   Assume  that  $$ \sup\{t>0: \psi(t) < \infty\} >1.  $$

  We shall consider the case when $(X_n)$ is null-recurrent and sub-diffusive.  Lyons and Pemantle  \cite{lyons-pemantle}   gave  a  precise recurrence/transience criterion for  randomly biased walks on an arbitrary infinite  tree. Their results,  together with   Menshikov and Petritis \cite{menshikov-petritis} and Faraud \cite{faraud},  imply that  $(X_n)$ is null recurrent if and  only if  $\inf_{0\le t\le  1} \psi(t)=  0$ and $\psi'(1)\le 0$.  There are two different situations in the null-recurrent case:  Either  $\psi'(1) = 0$, then $(X_n)$ has a slow-movement behavior  (i.e. $|X_n|$ is in the logarithmic scale, see \cite{gyzbiased} for the maximal displacement of $X$,  and see  \cite{yzlocaltree} for the localization of $X_n$ and the study of the local times processes), or $ \psi'(1)< 0$, then $(X_n)$ is sub-diffusive (in the sense of \eqref{polynomial}, see  \cite{yztree}).   We assume throughout  this paper 
   \begin{equation}\label{hyp1} \inf_{0\le t\le  1} \psi(t) =0  \qquad \hbox{and} \qquad  \psi'(1)< 0.\end{equation}

 \noindent Let us introduce a parameter \begin{equation*}   \kappa:= \inf\{t>1: \psi(t)=0\} \, \in (1, \infty], \end{equation*}
 
 \noindent with  $\inf\emptyset:=\infty$.  We furthermore assume the following conditions:  \begin{equation}\label{hyp2} 
 	\begin{cases}
  \E \Big( \sum_{i=1}^\nu \, A_i\Big)^\kappa + \E\Big( \sum_{i=1}^\nu A_i^\kappa \log_+ A_i\Big) < \infty ,  	\qquad & \hbox{\rm if $1 < \kappa \le 2$} \, ,
    \\
    \E \Big( \sum_{i=1}^\nu \, A_i\Big)^2 < \infty, \qquad   & \hbox{\rm if $ \kappa\in ( 2, \infty]$} \, , \end{cases}
	\end{equation}

 \noindent  with   $\log_+ x:= \max(0, \log x)$ for any $x>0$.


  \begin{figure}[h]
 \includegraphics[scale=1]{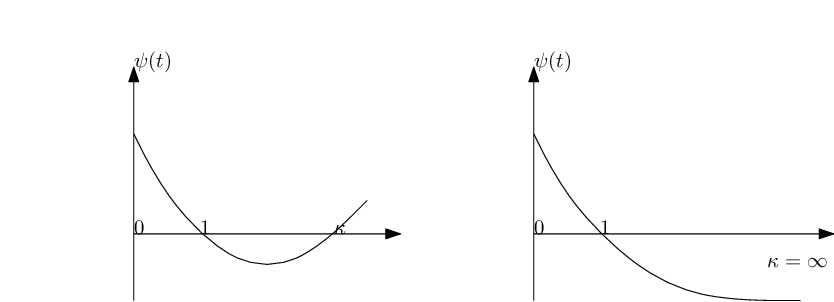} 
\caption{Case $\inf_{0\le t\le  1} \psi(t) =0$ and $\psi'(1)< 0$: $\kappa\in (1, \infty)$ and $\kappa=\infty$}
\end{figure}

 It was shown in  \cite{yztree} that if $\T$ is a regular tree  (i.e. $\nu$ equals some constant), then   \begin{equation}\label{polynomial} \lim_{n\to\infty} \frac1{\log n} \log \max_{0\le i\le n} |X_i| =  1 - \max( {1\over 2}, { 1 \over \kappa}), \qquad \hbox{$\p$-a.s.}.  \end{equation}

\noindent
When $\kappa$ is sufficiently large (say $\kappa \in (5, \infty]$), Faraud \cite{faraud} proved an invariance principle for the biased walk $X$, based on the techniques of Peres and Zeitouni \cite{peres-zeitouni}; some recent developments cover the whole region $\kappa \in (2, \infty]$ (see  A\"\i d\'ekon and de Raph\'elis \cite{elie-loic} for the convergence to Brownian tree).


The biased walk on a Galton-Watson tree has also attracted many attentions from other   directions: In the transient case,     A\"\i d\'ekon \cite{elie, elie2} dealt with a leafless Galton-Watson tree, whereas  Hammond  \cite{hammond} established  stable limit laws for the  walk on a supercritical Galton-Watson tree with leaves, which can be considered as a counterpart of Ben Arous, Gantert, Fribergh and Hammond  \cite{benarous-fribergh-gantert-hammond}. When the tree is sub-critical,  Ben Arous and Hammond \cite{benarous-hammond}  obtained power laws for the tails of $ E_\omega(T^+_\varnothing) $ and the convergence in law of $T^+_\varnothing$  under a suitable conditional probability, where $T^+_\varnothing$ denotes the return time to $\varnothing$: \begin{equation} \label{T+}
  T^+_\varnothing
  :=
  \inf\{n\ge1: X_n=\varnothing\}.
\end{equation}

\noindent  
In  the above-mentioned works \cite{hammond, benarous-fribergh-gantert-hammond, benarous-hammond},  the authors explored the  link between the biased walk $(X_n)$ and the {\it trap models} (cf.  Ben Arous and Cerny \cite{benarous-cerny}) to get various  scaling limits, and an important  step is the estimate on the return time to the trap entrance in their models.

We investigate here the return time $T^+_\varnothing$ in the scope of limit theorems for the local time process of the randomly biased walk $(X_n)$.    It turns out the parameter $\kappa$ plays a crucial role. Indeed, define $(M_n)$ by \begin{equation}\label{Mn}
	M_n
	:=
	\sum_{|x|=n} \prod_{\varnothing< y \le x} A(y), \qquad n\ge1, 
	\end{equation}

\noindent where here and in the sequel,  $|x|$ is the generation of $x$ in $\T$ and we adopt the partial order:  $y< x$ means that $y$ is ancestor of $x$ [we write $y\le x$ iff either $y< x$ or $y=x$].  Since $\psi(1)=0$, it is easy to check that $(M_n)$ is a martingale, which in the language of branching random walk is called the additive martingale (cf. Shi \cite{stf} further properties on $(M_n)$).   Define $$ \P^*(\bullet)
	:= 
	\P\Big(\bullet \, |\, \T=\infty\Big),$$

\noindent where $\{\T=\infty\}$ denotes  the event that the system survives forever.  Let $ M_\infty
	:= 
	\lim_{n \to \infty} M_n$ be the almost sure limit of the  nonnegative martingale  $(M_n)$. Then under \eqref{hyp1} and \eqref{hyp2} [the condition \eqref{hyp2} is more than necessary to ensure the non-triviality of $M_\infty$], $\P^*$-a.s.   $M_\infty >0$. 
Moreover, by  Liu \cite{liu00} Theorems 2.0 and 2.2,  when $\kappa\in (1, 2]$,     \begin{equation}
	 \P\Big( M_\infty > x \Big) \, \asymp\, x^{-\kappa},  \qquad x\to \infty, \label{cM}
\end{equation}
\noindent  where here and in the sequel, we denote by   $f(x) \asymp g(x)$  (resp: $f(x) \sim g(x)$) as $x\to x_0$ if  $0< \liminf_{x\to x_0}  f(x)/g(x) \le  \limsup_{x\to x_0}  f(x)/g(x) < \infty$ (resp: $\lim_{x\to x_0} f(x)/g(x)=1$)   [in fact \eqref{cM} holds for all $\kappa \in (1, \infty)$ under suitable integrability conditions, but here we only need it for the case $\kappa\in (1, 2]$].

The main estimate on the return time reads as follows.   

\begin{theorem}\label{t:main} Assume \eqref{hyp1}, \eqref{hyp2}.  We have $\P^*(d \omega)$-a.s., 

(i) if $1< \kappa<2$, then $$  P_{ \omega}\big(T^+_\varnothing >  n\big )  \, \asymp\, n^{-1/\kappa};$$

(ii) if $\kappa=2$, then  $$  P_{ \omega}\big(T^+_\varnothing >  n\big )  \, \asymp\, (n \log n)^{-1/2};$$

(iii) if $2 <  \kappa  \le  \infty $, then 
	$$
	\frac{1}{ \omega(\varnothing, {\buildrel \leftarrow \over \varnothing}) \, M_\infty}\, P_{ \omega}\big(T^+_\varnothing >  n\big )  \sim       
 c_1\, n^{-1/2},  	$$
	with  $$	c_1:=    \Big( \frac{2}{\pi}\,  \frac{1-\E(\sum_{i=1}^\nu A_i^2)}{\E(\sum_{1\le i\not= j\le \nu} A_i A_j)} \Big)^{1/2},$$

\noindent where ${\mathbb B}$ denotes  the Beta function. \end{theorem}

\medskip
As a consequence, we get  the asymptotic behaviors of  the local times process:  $$ L_n^x:= \sum_{i=1}^n 1_{(X_i=x)}, \qquad n\ge1, x\in \T.$$

\noindent  We shall restrict our attentions to the local times at the root.  It was implicitly contained  in \cite{yztree, andreoletti-debs1} that for any $ \kappa \in (1,  \infty]$, $\p$-almost surely on $\{\T=\infty\}$, $$ L_n^\varnothing= n^{\max(1/\kappa, 1/2)+o(1)}.$$

Based on Theorem \ref{t:main}, we can get more precise information on $L_n^\varnothing$. 

\begin{corollary}\label{c:main} Under the same assumptions as in Theorem \ref{t:main},  $\P^*(d\omega)$-a.s.,  we have that under $P_\omega$ as $n \to \infty$: 

(i) if $1< \kappa <2$, then $  \frac{L_n^\varnothing}{n^{1/\kappa}} $ is tight on $(0, \infty)$;
 

(ii) if $\kappa=2$, then $  \frac{L_n^\varnothing}{\sqrt{n \log n}}$ is tight on $(0, \infty)$;


(iii) if $ 2 <  \kappa  \le  \infty $, then $$ \frac{L_n^\varnothing}{\sqrt n}\, \wcv\,  \frac{1}{\omega(\varnothing, {\buildrel\leftarrow \over \varnothing}) M_\infty}\, \frac{2^{1/2}}{c_1\, \pi^{1/2}} \, |{\cal N}| ,$$

\noindent  where ${\cal N}$ denotes a standard gaussian random variable, centered and with variance $1$. 

\end{corollary}

\medskip
By the classical fluctuation theory on the random walk in the domain of attraction, it is straightforward  to deduce from Theorem \ref{t:main} the almost sure limits on $L_n^\varnothing$: for instance,  we have the following law  of iterated logarithm:

\begin{corollary}\label{c:lil}   Under the same assumptions as in Theorem \ref{t:main},   for any $\kappa \in (1, \infty]$,   $\P^*(d\omega)$-a.s.,  we have that under $P_\omega$-a.s.,       $$  	\limsup_{n \to \infty} \frac{L_n^\varnothing}{f_\kappa(n)}  \in (0, \infty),  $$

\noindent where $$ f_\kappa(n):=  \begin{cases}
 n^{1/\kappa} (\log\log n)^{1-1/\kappa},
	 \qquad & \hbox{\rm if $1 < \kappa < 2$} \, ,
    \\
    n^{1/2} (\log n)^{1/2}  (\log\log n)^{1/2},  \qquad & \hbox{\rm if $ \kappa= 2$}\,, 
    \\
 n^{1/2} ( \log\log n)^{1/2},  \qquad  & \hbox{\rm if $ \kappa\in (2, \infty]$} \, .
    \end{cases} $$
  \end{corollary}



Combining the estimates on the local times and the reversibility of the biased walk, we  obtain the following estimates on  the local probability.

\begin{corollary}\label{c:local} Under the same assumptions as in Theorem \ref{t:main},    $ \P^*(d\omega)$-a.s., for $n\to\infty$ along the sequence of even numbers,  

(i) if $1< \kappa <2$, then $$ P_\omega\Big( X_n= \varnothing\Big) \, \asymp \,  n^{-1+1/\kappa}; $$


(ii) if $\kappa=2$, then  $$ P_\omega\Big( X_n= \varnothing\Big)  \, \asymp \,  n^{-1/2} \, (\log n)^{1/2}; $$


(iii) if $ 2 <  \kappa  \le  \infty $, then  $$ P_\omega\Big( X_n= \varnothing\Big) \sim  \frac{1}{\omega(\varnothing, {\buildrel\leftarrow \over \varnothing}) M_\infty}  \,  \frac{2}{\pi c_1}\,   n^{-1/2} .$$

\end{corollary}

\section{Outline of the proof}

    For any $x \in \T$,  let  $P_{x, \omega}$ be the law of the biased walk $X$ starting from $X_0:=x$.  Denote by $E_{x, \omega}$ the expectation under the probability measure $P_{x, \omega}$.  In particular, we have $P_{\varnothing, \omega}\equiv P_\omega$ and $E_{\varnothing, \omega} \equiv E_\omega$.  Let  $$T_x:= \inf \{n\ge 0: X_n= x\}, \qquad x\in \T,$$ 

\noindent be  the first hitting time of $x$.    Clearly for $n >2$,   \begin{eqnarray}
	P_\omega\Big(T^+_\varnothing> n\Big)
	&=& 
	\sum_{|u|=1} \omega(\varnothing, u) P_{u, \omega}\Big( T_\varnothing> n-1\Big) \nonumber
	\\
	&=&
	\omega(\varnothing, {\buildrel \leftarrow \over \varnothing}) \sum_{|u|=1} A(u) P_{u, \omega}\Big( T_\varnothing> n-1\Big) .  \label{pt+1}
\end{eqnarray}

\noindent   By Tauberian theorems, the asymptotic behaviors of $P_{u, \omega}( T_\varnothing> n-1)$, are characterized by that of $ E_{u, \omega} \big( \ee^{- \lambda T_\varnothing}\big)$ as $\lambda \to 0$.  More generally, we define  for any $\lambda>0$ and $x \in \T$,  \begin{equation} \label{betax}	\beta_\lambda(x)
	:=   
	 1- E_{x, \omega} \Big( \ee^{- \lambda (1+ T_{{\buildrel\leftarrow \over x}})}  \Big), \qquad x \in \T, 
	 \end{equation}

\noindent  where as before, ${\buildrel\leftarrow \over x}$ denotes the parent of $x$. 
It is easy to see that $\beta_\lambda(\cdot)$ satisfies the following recursive iteration equations: 

\begin{fact} \label{f:beta} For any $x\in \T$ and $\lambda>0$, we have
	$$ \beta_\lambda(x)
	= 
	\frac{(1-\ee^{-2\lambda}) + \sum_{i=1}^{\nu_x} A(x^{(i)}) \beta_\lambda(x^{(i)})}{1+ \sum_{i=1}^{\nu_x} A(x^{(i)}) \beta_\lambda(x^{(i)})}.$$
\end{fact}

\medskip We mention  that conditioned on $\big((A(x^{(i)}))_{ 1\le i \le \nu_x}, \nu_x\big)$, $(\beta_\lambda(x^{(i)}), 1\le i \le \nu_x)$ are i.i.d. and are distributed as $\beta_\lambda(\varnothing)$. 

\medskip
{\noindent\bf Proof of Fact \ref{f:beta}.} This fact is an easy application of Markov property. We give the proof for the sake of completeness. For  use later,  we define  for any $n\ge1, \lambda>0$ and $x \in \T$ and $|x| \le n$,  \begin{equation} \label{betanx} \beta_{n, \lambda}(x)
	:= 
	1- E_{x, \omega} \Big( \ee^{- \lambda (1+ T_{{\buildrel\leftarrow \over x}})} 1_{(\tau_n > T_{{\buildrel\leftarrow \over x}})}\Big),
	\end{equation}

\noindent where  $$\tau_n:=\inf\{k\ge0: |X_k|=n\}, $$

\noindent denotes the first time when $X$ hits the $n$-th generation of the tree $\T$.

Clearly $\beta_{n, \lambda}(x)= 1$ for all $|x|=n$ and for $|x|< n$, we have by the Markov property that  \begin{eqnarray*}
	\beta_{n, \lambda}(x) &=& 1- \Big( \sum_{i=1}^{\nu_x} \omega(x, x^{(i)}) \ee^{-\lambda}\, E_{\omega, x^{(i)}} \ee^{-\lambda(1+T_{{\buildrel\leftarrow \over x}})} 1_{(T_{{\buildrel\leftarrow \over x}}< \tau_n)} + \omega(x, {\buildrel\leftarrow \over x}) \ee^{-2 \lambda}\Big) 
	\\ 
	&=& 1- \Big( \sum_{i=1}^{\nu_x} \omega(x, x^{(i)}) (1- \beta_{n, \lambda}(x^{(i)}) (1- \beta_{n, \lambda}(x)) + \omega(x, {\buildrel\leftarrow \over x}) \ee^{-2 \lambda}\Big) . \end{eqnarray*}

\noindent After simplifications, we get that \begin{equation} \label{beta-n}
	\beta_{n, \lambda}(x)= \frac{(1-\ee^{-2\lambda}) + \sum_{i=1}^{\nu_x} A(x^{(i)}) \beta_{n, \lambda}(x^{(i)})}{1+ \sum_{i=1}^{\nu_x} A(x^{(i)}) \beta_{n, \lambda}(x^{(i)})}, \qquad |x| < n.
\end{equation}

 \noindent Letting $n\to \infty$,  $\beta_\lambda(x)= \lim_{n \to\infty}\beta_{n,\lambda}(x)$ and we get Fact \ref{f:beta}. \hfill$\Box$

\medskip

For brevity, we make a change of variable  $ \varepsilon= 1- \ee^{-2\lambda},$ by defining 
 \begin{equation}\label{Bex}
  B_\varepsilon(x)
  :=
   \sum_{i=1}^{\nu_x} A(x^{(i)}) \beta_{\frac12 \log 1/(1-\varepsilon)}(x^{(i)}), \qquad x \in \T, \qquad 0< \varepsilon< 1,
   \end{equation}

\noindent then  \begin{equation} \label{Bx}
	B_\varepsilon(x)= \sum_{i=1}^{\nu_x} \, A(x^{(i)})\, \frac{\varepsilon  + B_\varepsilon(x^{(i)}) }{1+ B_\varepsilon(x^{(i)}) },
\end{equation}

\noindent where  as for $\beta_\lambda(x)$, conditioned on $\big((A(x^{(i)}))_{ 1\le i \le \nu_x}, \nu_x\big)$,  $(B_\varepsilon(x^{(i)}), 1\le i \le \nu_x)$ are i.i.d. and  are distributed as $B_\varepsilon(\varnothing) $.    In view of \eqref{betax}, we remark that \begin{equation}\label{Bepsilon}
B_\varepsilon(\varnothing)
=
\sum_{|u|=1}  A(u)\, \Big( 1- E_{u,  \omega}\big((1-\varepsilon)^{\frac{1+T_\varnothing}{2}}\big)\Big).
 \end{equation}
The main  estimate in the proof of Theorem \ref{t:main} will be the following result:

\begin{proposition}\label{p:B}   Assume  \eqref{hyp1}, \eqref{hyp2}.  As $\varepsilon \to 0$, the following  convergence holds $\P$-almost surely as well as in $L^p(\P)$ for any $1< p < \min (\kappa, 2)$: $$
\frac{B_\varepsilon(\varnothing)}{\E(B_\varepsilon(\varnothing))} \, \to\, M_\infty.$$

\noindent Moreover, as $\varepsilon\to0$, 

(i)  if $1 < \kappa< 2$, then $$ \E(B_\varepsilon(\varnothing)) \, \asymp \, \varepsilon^{1/\kappa};  $$

(ii) if $\kappa=2$, then   $$ \E(B_\varepsilon(\varnothing)) \, \asymp \, \Big(\frac\varepsilon{\log \frac1\varepsilon}\Big)^{1/2} ; $$

(iii)  if $2 < \kappa   \le \infty$, then $$  \E(B_\varepsilon(\varnothing)) \, \sim \, c_2\, \varepsilon^{1/2}, $$

\noindent where $c_2:=  \Big( \frac{1-\E(\sum_{i=1}^\nu A_i^2)}{\E(\sum_{1\le i\not= j\le \nu} A_i A_j) } \Big)^{1/2}$.  

\end{proposition}

\medskip
\noindent Recall that $\P$-a.s., $\{ M_\infty>0\}=\{\T=\infty\}$.  It is straightforward to see that on $\{T=\infty\}^c$, the biased walk $X$ is a Markov chain with finite states, hence $B_\varepsilon(\varnothing) =O(\varepsilon)$ as $\varepsilon\to0$.

\medskip

Let us give the proofs of Theorem \ref{t:main} and Corollaries \ref{c:main}, \ref{c:lil} and \ref{c:local}, by admitting Proposition \ref{p:B}:
\medskip

{\noindent\bf Proof of Theorem \ref{t:main}.} 
 By \eqref{betax} and \eqref{Bex}, we deduce from the usual Abel transform  that if $\lambda>0$ is such that $\varepsilon= 1-\ee^{-2\lambda}$, then $$ B_\varepsilon(\varnothing)
 	=
	(1-\ee^{-\lambda})\, \sum_{|u|=1} A(u) \sum_{k=0}^\infty \, \ee^{-\lambda k}\, P_{u, \omega}\big(T_\varnothing \ge k\big).$$

 In view of \eqref{pt+1}, $\sum_{|u|=1} A(u) P_{u, \omega}(T_\varnothing \ge k ) =  P_{ \omega}(T^+_\varnothing >  k )  / \omega(\varnothing, {\buildrel \leftarrow \over \varnothing})$. It follows that  \begin{equation}\label{laplace1}
 	\sum_{k=0}^\infty \, \ee^{-\lambda k}\, P_{ \omega}\big(T^+_\varnothing >  k\big )  = \omega(\varnothing, {\buildrel \leftarrow \over \varnothing}) \, \frac{B_\varepsilon(\varnothing)}{1-\ee^{-\lambda}}, \end{equation}
	
	\noindent with $\varepsilon= 1-\ee^{-2\lambda}$.

	When $2 < \kappa \le \infty$, we deduce from Proposition \ref{p:B} (iii) that $\P^*(d \omega)$-a.s., $B_\varepsilon(\varnothing) \sim c_2 M_\infty \varepsilon^{1/2}$, which according to  the Tauberian theorem (\cite{feller}, pp. 447, Theorem 5),  yields Theorem \ref{t:main} (iii).

	It remains to deal with the cases $\kappa \in (1, 2]$. For notational brevity, we define for any $0<\varepsilon\le 1$,  \begin{equation}\label{r} 
	r(\varepsilon):= 
	\begin{cases}
	\varepsilon^{1/\kappa}, \qquad & \mbox{if $\kappa \in (1, 2)$}, \\
	\sqrt{\varepsilon/\log(\ee/\varepsilon)} , \qquad & \mbox{if $\kappa =2$}.
	\end{cases}
	\end{equation}
	
	For any $n\ge1$, by  \eqref{laplace1} and with $\varepsilon= 1-\ee^{-2\lambda}$,  $$
	B_\varepsilon(\varnothing) 
	  \ge 
	  ( 1- \ee^{-\lambda}) \,  \sum_{k=0}^n \, \ee^{-\lambda k}\, P_{ \omega}\big(T^+_\varnothing >  n\big ) 
	   =  (1- \ee^{- \lambda (n+1)} )\,  P_{ \omega}\big(T^+_\varnothing >  n\big ) .$$

	     Taking $\lambda= 1/(n+1)$ (then $\varepsilon= 1- \ee^{-2/(n+1)}$), we deduce from Proposition \ref{p:B}  that  \begin{equation}\label{uppT>n} P_{ \omega}\big(T^+_\varnothing >  n\big ) \le  \frac1{1- \ee^{-1}} \, B_\varepsilon(\varnothing) \le  c_3 \, r(1/n) , \qquad \forall \, n \ge 1,\end{equation}

\noindent where $c_3\equiv c_3(\omega) \in (0, \infty)$ only depends on the environment $\omega$. 

To get the lower bound for $P_{ \omega}\big(T^+_\varnothing >  n\big ) $, we use  \eqref{uppT>n} in the left-hand-side of \eqref{laplace1} and  obtain  that  for any $n\ge1$, 
\begin{eqnarray*}
 \sum_{k=0}^\infty \, \ee^{-\lambda k}\, P_{ \omega}\big(T^+_\varnothing >  k\big ) 
&\le &
1+ \sum_{k=1}^n \, \ee^{-\lambda k}\,  c_3\,  r(1/k) +  \sum_{k=n+1}^\infty \, \ee^{-\lambda k}\,   P_{ \omega}\big(T^+_\varnothing >  n\big ) 
\\
&\le&
1+ \sum_{k=1}^n \, \ee^{-\lambda k}\,  c_3\,  r(1/k) +  \frac{1}{1- \ee^{-\lambda}} P_{ \omega}\big(T^+_\varnothing >  n\big ). 
\end{eqnarray*}

Take $\lambda= \delta/n$ with a small $\delta>0$ [the value of $\delta$ will be determined later]. Elementary computations yield  the existence of some positive constant $c_4(\delta)$ satisfying that $c_4(\delta) \to 0$ as $\delta\to 0$ and such that for all $n\ge1$ and $\lambda= \delta/n$,  $$1+ \sum_{k=1}^n \, \ee^{-\lambda k}\,   r(1/k)   \le   c_4(\delta)\, r(\lambda)/\lambda.$$ 

By \eqref{laplace1} and Proposition \ref{p:B}, $\P^*(d \omega)$-a.s. there exists some $c_5\equiv
c_5(\omega) \in (0, \infty)$ such that for any $0< \lambda < 1$, 
 $\sum_{k=0}^\infty \, \ee^{-\lambda k}\, P_{ \omega}\big(T^+_\varnothing >  k\big )  \ge c_5\,  r(\lambda)/\lambda$.

  Choose (and then fix) $\delta$ sufficiently small such that $c_4(\delta) \le \frac12 \, c_5/c_3$, we have    that for any $n\ge 1$ and   $\lambda= \delta/n$,  \begin{equation}\label{lowT>n} P_{ \omega}\big(T^+_\varnothing >  n \big ) 
 \ge \frac12 c_5 \,  r(\lambda)  \frac{1- \ee^{-\lambda}}{\lambda}  \ge c_6\, r(1/n), \end{equation}
 
 \noindent for some $c_6 \equiv c_6(\omega, \delta) \in (0, \infty)$. Then \eqref{lowT>n}  together with \eqref{uppT>n} prove  Theorem \ref{t:main} (i) and (ii). 	 \hfill$\Box$

 \medskip
 {\noindent\bf Proofs of Corollaries  \ref{c:main} and   \ref{c:lil}.}  Fix a realization of environment $\omega$ such that  \eqref{uppT>n},   \eqref{lowT>n} and Theorem \ref{t:main} (iii)   hold.

   Define for $k\ge1$, $$ T^{(k)}_\varnothing:= \inf\{ n> T^{(k-1)}_\varnothing : X_n = \varnothing\}, $$
 
 \noindent   the $k$-th return to $\varnothing$ (with $T^{(0)}_\varnothing:=0$).  Then for any $n\ge1$, $L_n^\varnothing= \inf\{k\ge1: T^{(k)}_\varnothing \ge  n\}$. 
 
    Under $P_\omega$,  $T^{(k)}_\varnothing$ is the sum of $k$ i.i.d. copies of $T^+_\varnothing$, which  is in the domain of attraction of a normal law when $2< \kappa \le \infty$. For the cases $\kappa \in (1, 2]$, we use the stochastic dominance 
    between $T^+_\varnothing$ and a stable r.v.  of index $\max(1/\kappa, 1/2)$ as follows:  Let $\eta$ and $\widehat \eta$ be r.v. taking values in positive  integers such that for any $n\ge1$, $$ P_\omega\big( \eta>n)= \min(1, c_3 r(1/n)), \qquad P_\omega\big( \widehat \eta>n)= c_6 r(1/n),$$
    
    \noindent  where $r(\cdot)$, $c_3$ and $c_6$ are given in \eqref{r},  \eqref{uppT>n},   \eqref{lowT>n}  respectively, and by \eqref{lowT>n},  $c_6 r(1/n)\le 1$.  Then $\eta$ and $\widehat\eta$ are in the domain of attraction of a stable law of index $\max(1/\kappa, 1/2)$.   By \eqref{uppT>n} and  \eqref{lowT>n}, we may construct $\eta$,  $\widehat \eta$ and $T^+_\varnothing$ in a same   probability space  such that $P_\omega$-a.s., $\eta\ge T^+_\varnothing\ge \widehat \eta$.  
    
     Then for the cases $\kappa \in (1, 2]$,  in an eventually enlarged probability space we may construct two random walks $(\Xi_k)_{k\ge1}$  and  $(\widehat\Xi_k)_{k\ge1}$ such that for any $k\ge 1$,   $\Xi_k$ (resp: $\widehat \Xi_k$) is the sum of $k$ i.i.d. copies of $\eta$ (resp: $\widehat\eta$), and    $P_\omega$-a.s., $\Xi_k\ge T^{(k)}_\varnothing\ge \widehat \Xi_k. $ It follows that   $P_\omega$-a.s.,  $$\Xi^{-1}_n \le L_n^\varnothing \le \widehat\Xi^{-1}_n,  \qquad  \forall n\ge1,$$ with $\Xi^{-1}_n :=\inf\{k\ge1: \eta_k\ge n\}$ and similar definition   for $\widehat\Xi^{-1}_n$. 
    
    When $\kappa\in (2, \infty]$,   we take  $\Xi_k= T^{(k)}_\varnothing$ (then $\Xi_n^{-1}= L_n^\varnothing$ in this case) for any $k\ge1$.  Therefore for all cases $\kappa \in (1, \infty]$,   it suffices to prove that the conclusions of Corollaries  \ref{c:main} and   \ref{c:lil} hold for the process $(\Xi^{-1}_n)_{n\ge1}$ in lieu of $(L_n^\varnothing)_{n\ge 1}$. We may assume  in the sequel that  for some $0<\alpha<1$ and  some slowly varying function $\ell(n)$,\footnote{If $1 < \kappa\le 2$, $\alpha=1/\kappa$ and $\ell(n)= c_3 \, \Gamma(1-1/\kappa) n^{1/\kappa} r(1/n)$, whereas if $\kappa \in (2, \infty]$, $\alpha=1/2$ and  $\ell(n)= c_1 \Gamma(1/2)  \omega(\varnothing, {\buildrel\leftarrow \over \varnothing}) M_\infty$.}    \begin{equation} \label{alpha1} P_\omega(\Xi_1 > n) 
 	\sim 
	\frac1{\Gamma(1-\alpha)}\, n^{-\alpha} \ell(n) . \end{equation}
	
By  \cite{feller} (Theorem 2, pp.448), we have  that under $P_\omega$, $$ \frac{\Xi_k}{\big(k  \, \ell(k^{1/\alpha}) \big)^{1/\alpha}}\, \wcv\, {\cal S_\alpha},$$

\noindent with ${\cal S}_\alpha$ a positive stable variable of index $\alpha$ whose Laplace transform is given by $\e \ee^{-\lambda {\cal S}_\alpha}= \ee^{-\lambda^\alpha}$ for any $\lambda>0$.   In particular, ${\cal S}_{1/2}\law  \frac{1}{2 {\cal N}^2}$.

Applying  Fristed and Pruitt (\cite{FP71}, Theorem 5) to  the random walk $(\Xi_k)_{k\ge 1}$, we see that under $P_\omega$, $$ \limsup_{n\to\infty} \frac{\Xi^{-1}_n}{f_\alpha(n)} \in (0, \infty),$$

\noindent where $f_\alpha(\cdot)$ is given in   Corollary \ref{c:lil}.   This implies Corollary \ref{c:lil}.

Now  using the fact that $P_\omega( \Xi^{-1}_n \ge k) = P_\omega(\Xi_k\le n)$ for $k\ge1, n\ge 1$, we   deduce that for any $z>0$,
	$$  P_\omega\Big( \frac{ \Xi^{-1}_n}{n^\alpha/\ell(n)} \ge z\Big) 
	\to 
	\p\Big({\cal S}_\alpha \le z^{-1/\alpha}\Big).$$

\noindent It follows that under $P_\omega$,  \begin{equation}\label{wcvalpha}
	n^{-\alpha}\ell(n) \, \Xi^{-1}_n  \wcv ({\cal S}_\alpha)^{-\alpha},
 \end{equation}
 
 \noindent  which implies  Corollary \ref{c:main}. 
 \hfill$\Box$

\medskip

{\noindent\it Proof of Corollary \ref{c:local}.}  Under the framework   \eqref{alpha1}, we remark that  $n^{-\alpha}\ell(n) \, \Xi^{-1}_n $ is bounded in $L^p(P_\omega)$ for any $p>0$. In fact,   
$$E_\omega \Big( \Xi^{-1}_n \Big)^p \le \sum_{k=0}^\infty p \,  \, k^{p-1}\, P_\omega\Big(\Xi^{-1}_n \ge k\Big)  =\sum_{k=0}^\infty p  \, k^{p-1}\, P_\omega\Big( \Xi_k \le n \Big) .$$

\noindent Observe that  $P_\omega\big(  \Xi_k \le n \big)\le P_\omega\big( \Xi_1 \le n \big)^k  \le  \ee^{- k\, P_\omega ( \Xi_1 > n  )}$.  Hence $$E_\omega \Big( \Xi^{-1}_n \Big)^p \le \sum_{k=0}^\infty p  \, k^{p-1}\,  \ee^{- k\, P_\omega ( \Xi_1 > n  )}.$$

 \noindent Since $\sum_{k=0}^\infty p  \, k^{p-1}\, \ee^{- k x}\le c_7\,  x^{-p} $ for all $0< x\le1$ and some constant $c_7=c_7(p)>0$, we get that  $   P_\omega\big( \Xi_1> n \big) \times \Xi^{-1}_n $ is bounded in $L^p$ for any $p>0$.   This together with \eqref{wcvalpha} imply  that  \begin{equation}\label{EXin} 	E_\omega(  \Xi^{-1}_n)  \, 
 	\sim \,
	\e\big(({\cal S}_\alpha)^{-\alpha}\big)\, \frac{n^\alpha}{\ell(n)}, \qquad n\to \infty. \end{equation}
 
 Under $P_\omega$,  the Markov chain $X$ is reversible and  it is well-known (see e.g. Saloff-Coste (\cite{Saloff97}, Lemma 1.3.3 (1), page 323)) that $k \to P_\omega(X_{2k}=\varnothing)$ is non-increasing.

When $\kappa \in (2, \infty]$, $ \Xi^{-1}_n=L_n^\varnothing$ and  \eqref{EXin} says that $$ E_\omega( L_{2n}^\varnothing)  \, 
 	\sim \,
	 \frac{1}{\omega(\varnothing, {\buildrel\leftarrow \over \varnothing}) M_\infty}  \,  \frac{2}{\pi c_1}\,   (2n)^{1/2} , \qquad n\to \infty.$$ 
 
Since $E_\omega( L_{2n}^\varnothing) = \sum_{k=1}^n P_\omega(X_{2k}=\varnothing)$,   the Tauberian theorem (\cite{feller}, formula (5.26), pp.447)  yields  Corollary \ref{c:local}  (iii).

It remains to treat the cases $\kappa\in (1, 2]$.   Recalling that  $P_\omega$-a.s., $\Xi^{-1}_n \le L_n^\varnothing \le \widehat\Xi^{-1}_n,     \forall n\ge1$.  By \eqref{EXin}  and its analogue for $\widehat \Xi_n^{-1}$,  there are $c_8\equiv c_8(\omega) \in (0, \infty)$ and $c_9\equiv c_9(\omega) \in (0, \infty)$, such that for all $n \ge 1$, $$  \frac{c_8}{r(1/n)} \le E_\omega (L_{2n}^\varnothing) \le   \frac{c_9}{r(1/n)}.$$

Note that $E_\omega( L_{2n}^\varnothing) = \sum_{k=1}^n P_\omega(X_{2k}=\varnothing) \ge n P_\omega(X_{2n}=\varnothing) $, we get that  \begin{equation}\label{uppreturnproba} P_\omega(X_{2n}=\varnothing) \le \frac{c_9}{n \,r(1/n)}.\end{equation}

For the lower bound, we choose and fix a sufficiently small $\delta\equiv \delta(\omega)>0$ such that for any $n\ge1$, $r(1/(\delta n)) \ge  \frac{2c_9}{c_8} r(1/n)$.\footnote{Strictly speaking we use  $\lfloor \delta n\rfloor$  the integer part of $\delta n$ in this whole paragraph.}   Then for all $n\ge1$, $E_\omega(L_{2n}^\varnothing)- E_\omega(L_{2\delta n}^\omega) \ge  \frac{c_8}{r(1/n)} - \frac{c_9}{r(1/(\delta n))}  \ge  \frac{c_8}{2r(1/n)} $.  By the monotonicity, $E_\omega(L_{2n}^\varnothing)- E_\omega(L_{2\delta n}^\omega)  =\sum_{k=\delta n +1}^n P_\omega(X_{2k}=\varnothing) \le n (1-\delta) P_\omega(X_{2 \delta n }=\varnothing)$.    It follows that for all $n\ge1$, $P_\omega(X_{2 \delta n }=\varnothing) \ge \frac{c_8}{2 \, n\, r(1/n)} $ which together with \eqref{uppreturnproba} imply Corollary \ref{c:local} (i) and (ii).   \hfill$\Box$

\medskip

The rest of this paper is devoted to the proof of  Proposition \ref{p:B}, which will be mainly driven by the recursive equations \eqref{Bx}.    Aldous and Bandyopadhyay \cite{AldousB} pointed out the variety of contexts where the recursive equations have arisen in various models on tree, see also Peres and Pemantle \cite{pemantle-peres2} for the studies of  a family of concave recursive iterations using the potential theory.    We analyze  here the equations \eqref{Bx} in the spirit of \cite{yztree} by establishing some comparison inequalities on the concave iteration.


The key point in the proof of Proposition \ref{p:B}  will be the asymptotic behavior of $\E(B_\varepsilon(\varnothing))$.  In Section \ref{s:recu},  we   obtain the lower bound for $\E(B_\varepsilon(\varnothing))$ for all $\kappa \in (1, \infty]$ and get the convergence in law for $\varepsilon^{-1/2} B_\varepsilon(\varnothing)$ for  $\kappa\in (2, \infty]$. The upper bound of  $\E(B_\varepsilon(\varnothing))$ will be presented  in  Section \ref{s:4}, where we shall complete the proof of   Proposition \ref{p:B}  by establishing the almost sure convergence of $\frac{B_\varepsilon(\varnothing)}{\E(B_\varepsilon(\varnothing))}$ to $M_\infty$.

Throughout the rest of this paper, $(c_i)_{10 \le i \le 23}$ denote some positive constants whose values may  depend on  some parameters such as $\kappa$ and  $p \in (1, \kappa)$. 

\section{Concave recursions  on trees}\label{s:recu}

Let $0< \varepsilon< 1$. By \eqref{Bx}, $B_\varepsilon(\varnothing)$ is a nonnegative solution of the following equation in law:  \begin{equation}\label{Blaw}  
	B_\varepsilon 
	\law 
	\sum_{i=1}^\nu A_i {\varepsilon + B_\varepsilon(i)\over 1+ B_\varepsilon(i)}, \end{equation}

\noindent  where as before, $(A_i, 1\le i \le \nu)\equiv  (A(x), |x|=1)$ and conditioned on $(A_i, 1\le i \le \nu)$,  $B_\varepsilon(i)$ are i.i.d., and are distributed as $B_\varepsilon$. We recall that   $\E \big( \sum_{i=1}^\nu A_i\big)=1$ and $\E \big( \sum_{i=1}^\nu A_i^\kappa\big)=1$ if $\kappa < \infty$.

It is easy to get the uniqueness among the nonnegative solutions. Indeed,  If $B_\varepsilon$ and $\widetilde B_\varepsilon$ are two nonnegative solutions, then in some enlarged probability space, we can find a coupling of  $(A_i, 1\le i \le \nu)$, $(B_\varepsilon,  B_\varepsilon(i), 1\le i \le \nu)$ and  $(\widetilde B_\varepsilon,  \widetilde B_\varepsilon(i), 1\le i \le \nu)$ such that the equation \eqref{Blaw} hold a.s. for $B_\varepsilon$ and $\widetilde B_\varepsilon$. Since  $B_\varepsilon$ is stochastically dominated by $\sum_{i=1}^\nu A_i$ hence integrable, we get that   $\E | B_\varepsilon- \widetilde B_\varepsilon| \le \E | \frac{\varepsilon+B_\varepsilon}{1+B_\varepsilon}- \frac{\varepsilon+\widetilde B_\varepsilon}{1+\widetilde B_\varepsilon}| \le (1-\varepsilon) \E | B_\varepsilon- \widetilde B_\varepsilon|$ which implies that $B_\varepsilon= \widetilde B_\varepsilon$ and the claimed uniqueness in law.  Therefore we write indistinguishably $B_\varepsilon\equiv B_\varepsilon(\varnothing)$.


 
This section is devoted to the asymptotic behaviors of $\E(B_\varepsilon)$ as $\varepsilon \to 0$.    Specifically, if $\kappa \in (2, \infty]$ which is the easier case, we shall obtain an exact asymptotic of $\E(B_\varepsilon)$ as $\varepsilon \to 0$, whereas for $\kappa \in (1, 2]$ we shall  get a lower bound, the corresponding  upper bound will be proved in Section \ref{s:4}.

First we check that $B_\varepsilon\to0$ in $L^1(\P)$. Notice  that $ \E (B_\varepsilon)= \E \frac{\varepsilon+B_\varepsilon}{1+B_\varepsilon}$ (since $\E \big( \sum_{i=1}^\nu A_i\big)=1$),  which after simplification gives that  \begin{equation}\label{B2}
	\E \Big( \frac{B_\varepsilon^2}{1+B_\varepsilon}\Big)
	= 
	\varepsilon\, \E\Big( \frac1{1+B_\varepsilon}\Big).\end{equation}

\noindent Therefore   $$ \E (B_\varepsilon)= \E \frac{\varepsilon+B_\varepsilon}{1+B_\varepsilon} \le \varepsilon + \E \frac{B_\varepsilon}{1+B_\varepsilon}  \le \varepsilon + \Big( \E \frac{B_\varepsilon^2}{(1+B_\varepsilon)^2} \Big)^{1/2} , $$

\noindent which in view of \eqref{B2} yield that  for any $\kappa \in (1, \infty]$,  \begin{equation}\label{EB1/2}
 \E (B_\varepsilon) \le  2  \varepsilon^{1/2}, \qquad 0<\varepsilon\le1.\end{equation}

The above upper bound is sharp (up to a constant) only in the case $\kappa \in (2, \infty]$. To obtain the lower bound on $\E(B_\varepsilon)$, we shall need  some inequalities on  the concave iteration.   Let us adopt the following notation in the rest of this paper:   $$\langle \xi \rangle:= \frac{\xi}{\E(\xi)},$$

\noindent for any nonnegative random variable $\xi$ with finite mean [as such, $\E \langle \xi\rangle^p= \frac{\E(\xi^p)}{(\E \xi)^p}$].

\begin{lemma}\label{L:inequality} Let $\phi: \r_+ \to \r_+$ be a convex $C^1$-function.  Let $\xi$ be an    nonnegative random variable $\xi$ with finite and positive mean. Suppose that there exists some $\delta>0$ such that $\E \phi( (1+\delta) \langle \xi\rangle)< \infty$. Then for  any $0\le \varepsilon < 1$, we have $$ \E  \phi\Big( \Big\langle \frac{\varepsilon+ \xi}{1+\xi} \Big\rangle  \Big) \le \E \phi\big(\langle\xi\rangle \big).$$
\end{lemma}

{\noindent\bf Proof:}  We shall use several times  the following inequality in \cite{yztree}, formula (3.3):     Let $x_0 \in \r_+$ and let  $I \subset \r _+$ be an open interval containing  $x_0$.  Assume that  $h: I \times \r_+ \to (0, \infty)$ is a Borel function such that  ${\partial h\over \partial x}$ exists and 
\begin{itemize}

 \item $\E [ h(x_0, \xi) ] <\infty$ and $\E \big[
       \phi \big( \big\langle h(x_0,\xi) \big\rangle
       \big)  \big] < \infty$;

 \item $\E[\sup_{x\in I} \big( | {\partial h\over
       \partial x} (x, \xi)| + |\phi' \big(\big\langle h(x,
       \xi) \big\rangle \big) | \,
       ({| {\partial h\over \partial x} (x, \xi) |
       \over \E ( h(x, \xi)) } + {h(x, \xi) \over
       [\E \langle h(x, \xi)\rangle]^2 } | \E ( {\partial
       h\over \partial x} (x, \xi) ) | )\big) ] <
       \infty$;

 \item both $y \to h(x_0, y)$ and $y \to {
       \partial \over \partial x} \log
       h(x,y)|_{x=x_0}$ are monotone on $\r_+$.

\end{itemize}

\noindent Then  depending on whether $h(x_0, \cdot)$ and ${\partial \over \partial x} \log h(x_0,\cdot)$ have the same monotonicity,
\begin{equation}
    \frac{d}{d x}\, \E  \phi\left(\big\langle
    h(x,\xi) \big\rangle\right)  
    \big|_{x=x_0} \ge 0, \qquad \hbox{\rm or}\qquad
    \le 0.
    \label{monotonie}
\end{equation}

Applying \eqref{monotonie} to $h(x,y):= \frac{x+y}{1+y}$, $0< x < 1$ and $y\ge0$. For any fixed $x_0\in (0, 1)$, $h(x_0, \cdot)$ is non-decreasing whereas $\frac{\partial}{\partial x} \log h(x_0, \cdot)= \frac{1}{x_0+\cdot}$ is non-increasing.  Therefore $x_0\in (0, 1) \mapsto \E\phi\big( \big\langle h(x_0, X)\big\rangle\big)$ is non-increasing. It follows that for any $0< \varepsilon< 1$, $$ \E \phi\Big( \Big\langle \frac{\varepsilon+ \xi}{1+\xi} \Big\rangle  \Big) \le \E \phi\Big( \Big\langle \frac{\xi}{1+\xi} \Big\rangle  \Big).$$

Now we take $h(x,y):= \frac{y}{1+xy}$ for $x\in (0,1)$ and $y\ge0$ in  \eqref{monotonie} and get that $x\in (0,1) \mapsto \E\phi\big( \big\langle \frac{\xi}{1+ x\, \xi} \big\rangle  \big) $ is non-increasing. Hence $\E \phi\big( \big\langle \frac{\xi}{1+\xi} \big\rangle  \big) \le \lim_{x\to 0} \E\phi\big( \big\langle \frac{\xi}{1+ x\, \xi} \big\rangle  \big) $.  For sufficiently small $x>0$, $\E( \frac{\xi}{1+ x\, \xi} ) > \E(\xi)/(1+\delta)$ which implies that $\big\langle \frac{\xi}{1+ x\, \xi} \big\rangle  \le (1+\delta) \big\langle  \xi \big\rangle $.  The convexity of $\phi$ implies that $\phi(\big\langle \frac{\xi}{1+ x\, \xi} \big\rangle  ) \le \max(\phi(0), \phi((1+\delta) \big\langle  \xi \big\rangle))$ which yields, by the dominated convergence theorem, that $ \lim_{x\to 0} \E\phi\big( \big\langle \frac{\xi}{1+ x\, \xi} \big\rangle  \big)  =\E \phi\big( \big\langle \xi \big\rangle  \big)$ and proves   the Lemma.  \hfill $\Box$

\begin{lemma} \label{l:Bp}  Assume  \eqref{hyp1} and \eqref{hyp2}. For any $p\in (1, 2] \cap (1,  \kappa)$, there exists some positive constant $c_{10}=c_{10}(p, \kappa) $ such that for any $0< \varepsilon < 1$,  \begin{eqnarray}  \E \Big( \big\langle B_\varepsilon\big\rangle^p\Big)  
	&\le & c_{10}. \,  \label{Bp1}
 \end{eqnarray}
\end{lemma}

{\noindent\bf Proof:}  The proof of  \eqref{Bp1} was already given in (\cite{yztree}, Proposition 5.1)  in the case that $\nu$ equals some integer larger than $2$.   The same proof can be adopted to the case of random $\nu$ and  
 we include the   proof here for the sake of completeness. Since $\E \big( \sum_{i=1}^\nu A_i\big)=1$, we have by the independence between $(B_\varepsilon(i))$ and $(A_i)$ that $$ \E (B_\varepsilon)= \E \Big( \frac{\varepsilon+B_\varepsilon}{1+B_\varepsilon}\Big),$$

 \noindent  which yields that  $$   \langle B_\varepsilon\rangle^p   =   \Big( \sum_{i=1}^\nu A_i  \Big\langle \frac{\varepsilon + B_\varepsilon(i)}{1+ B_\varepsilon(i)} \Big\rangle \Big)^p. $$

 \noindent We recall the following inequality due to Neveu \cite{neveu87}:  Let $k\ge1$ and let $\xi_1$, $\cdots$, $\xi_k$ be independent  non-negative random variables such that $\E(\xi_i^p)<\infty$ $(1\le i\le
  k)$. Then $$
  \E \left( \xi_1 + \cdots + \xi_k \right)^p  \le
  \sum_{i=1}^k \E(\xi_i^p) + \left(
  \sum_{i=1}^k \E \xi_i \right)^p.
  $$

 \noindent It follows that  
  \begin{eqnarray*}
 	\E  \Big[ \langle B_\varepsilon\rangle^p\Big] & \le& \E \sum_{i=1}^\nu A_i^p \, \Big\langle \frac{\varepsilon+B_\varepsilon}{1+B_\varepsilon}\Big\rangle^p + \E \Big(\sum_{i=1}^\nu A_i\Big)^p 
 	\\ &=& a_p \, \E \Big[ \Big\langle \frac{\varepsilon+B_\varepsilon}{1+B_\varepsilon}\Big\rangle^p\Big] + c_{11},
 \end{eqnarray*}

 \noindent where $c_{11}:= \E \big(\sum_{i=1}^\nu A_i\big)^p < \infty$ by the assumption \eqref{hyp2} and  \begin{equation}\label{ap} a_p :=  \E \sum_{i=1}^\nu A_i^p < 1, \end{equation}
 
 
 \noindent by the definition of $\kappa$. 
  Applying     Lemma \ref{L:inequality}  to $\phi(x):= x^p$, we get that $ \E  \Big\langle \frac{\varepsilon+B_\varepsilon}{1+B_\varepsilon}\Big\rangle^p \le  \E  \langle B_\varepsilon\rangle^p$, thus $\E  \langle B_\varepsilon\rangle^p \le \frac{c_{11}}{1-a_p}$, proving  \eqref{Bp1}.
\hfill $\Box$

\medskip

To get a lower bound of $\E(B_\varepsilon)$, we shall use the following comparison lemma:

\begin{lemma}\label{l:phia}  Assume \eqref{hyp1} and \eqref{hyp2}.   Let $p\in [1, \kappa) \cap (1, 2]$.  For any $a>0$, let $\phi_a(x):=   \big( \frac{x^2}{a+x} \big)^p $ for any $x\ge 0$.  We have $$
	\E  \phi_a\Big( \big\langle  B_\varepsilon \big\rangle  \Big)
	 \le  
	  \E \phi_a(M_\infty).$$
\end{lemma}

\medskip
{\noindent\bf Proof:}   It is elementary to check that  the function $\phi_a $ is  convex.   Moreover, for any $b\ge 0$ and $t >0$, the function $x \mapsto \phi_a(b +t x)$ is still convex. By   Lemma \ref{L:inequality}, we get that  for any $b\ge 0$ and $t >0$,  \begin{equation}\label{b+t} \E  \phi_a\left(b+ t \, \Big\langle  \frac{\varepsilon+ \xi }{1+\xi } \Big\rangle  \right) 
	\le
	 \E \phi_a(b + t\, \langle\xi \rangle).
	  \end{equation}


\noindent
 Recall \eqref{betanx}.  Choose $\lambda$ such that $1- \ee^{-2 \lambda}= \varepsilon$. Define $$B_{\varepsilon, n}(x)
	:=
	 \sum_{i=1}^{\nu_x} A(x^{(i)})  \beta_{n, \lambda}(x^{(i)}), \qquad \forall\,  |x|\le n.
	 $$

\noindent
Then $B_\varepsilon=B_\varepsilon(\varnothing) = \lim_{n\to \infty} B_{\varepsilon, n}(\varnothing)$,  $\P$-almost surely. For any $|x| < n$, we deduce from  \eqref{beta-n} that  $$ B_{\varepsilon, n}(x)= \sum_{i=1}^{\nu_x} A(x^{(i)}) \frac{\varepsilon+ B_{\varepsilon, n}(x^{(i)})}{1+B_{\varepsilon, n}(x^{(i)})} .$$

\noindent
Since $\E(\sum_{i=1}^{\nu_x} A(x^{(i)})) =1$,  we get that   $$ \langle B_{\varepsilon, n}(x)\rangle= \sum_{i=1}^{\nu_x}  A(x^{(i)})  \Big\langle\frac{\varepsilon+ B_{\varepsilon, n}(x^{(i)})}{1+B_{\varepsilon, n}(x^{(i)})}\Big\rangle .$$

\noindent
Applying \eqref{b+t} to $\xi= B_{\varepsilon, n}(x^{(1)})$, $t=A(x^{(1)})$ and  $b:= 1_{(\nu_x\ge2)}\,\sum_{i=2}^{\nu_x} A(x^{(i)})  \Big\langle\frac{\varepsilon+ B_{\varepsilon, n}(x^{(i)})}{1+B_{\varepsilon, n}(\varnothing^{(i)})}\Big\rangle$ and conditioning on $(t, b)$, we have  that $$ \E \phi_a(\langle B_{\varepsilon, n}(x)\rangle) \le \E \phi_a\Big( 1_{(\nu_x\ge2)}\,\sum_{i=2}^{\nu_x} A(x^{(i)})  \Big\langle\frac{\varepsilon+ B_{\varepsilon, n}(x^{(i)})}{1+B_{\varepsilon, n}(x^{(i)})}\Big\rangle + A(x^{(1)}) \langle B_{\varepsilon, n}(x^{(1)})\rangle\Big).$$

\noindent In the right-hand-side of the above inequality, applying  \eqref{b+t}  successively to $B_{\varepsilon, n}(x^{(2)}), ..., B_{\varepsilon, n}(x^{(\nu_x)})$  with obvious  choices of $t$ and $b$,   we get that for any $|x| < n$,  $$ \E \phi_a\Big(\langle B_{\varepsilon, n}(x)\rangle \Big) 
	\le 
	\E \phi_a\Big(  \sum_{i=1}^{\nu_x} A(x^{(i)})  \langle B_{\varepsilon, n}(x^{(i)})\rangle\Big).$$

\noindent
Notice by definition $B_{\varepsilon, n}(x)=  \sum_{i=1}^{\nu_x} A(x^{(i)})$ for $|x|=n-1$.  By iterating the above inequalities, we get that $$ \E \phi_a\Big(\langle B_{\varepsilon, n}(\varnothing)\rangle \Big) 
	\le 
	\E \phi_a\Big(  \sum_{|x|=n} \prod_{\varnothing< y\le x} A(y)  \Big) = \E \phi_a(M_n).$$

\noindent   Lemma  \ref{l:phia}  follows by letting $n \to \infty$.   \hfill$\Box$


\begin{lemma}\label{l:low}  Assume \eqref{hyp1}, \eqref{hyp2}.  We have \begin{eqnarray}
	 \liminf_{\varepsilon\to0} \,  \varepsilon^{-1/\kappa }\, \E(B_\varepsilon) \,
	 & >& \, 0 , \qquad  \hbox{if } 1< \kappa <2,  \label{kappa<2}
	 \\
	  \liminf_{\varepsilon\to0} \,  \Big(\frac{\log 1/\varepsilon}{\varepsilon}\Big)^{1/2} \, \E(B_\varepsilon) \, &>& \, 0, \qquad \hbox{if }   \kappa =2, \label{kappa=2}
	  \\
	    \liminf_{\varepsilon\to0} \,  \varepsilon^{-1/2} \, \E(B_\varepsilon) \, &\ge & \, c_2, \qquad \hbox{if }   \kappa >2, \label{kappa>2}
\end{eqnarray} 
where $c_2:=\big( \frac{1-\E(\sum_{i=1}^\nu A_i^2)}{\E(\sum_{1\le i\not= j\le \nu} A_i A_j)} \big)^{1/2}$ as in Proposition  \ref{p:B}.  
\end{lemma}

\medskip
{\noindent\bf Proof:}    If $\kappa\in (2, \infty]$, we remark that  $\E(\sum_{i=1}^\nu A_i^2) < 1$ and    $\E(M_\infty^2)
	=
	 \frac{\E(\sum_{1\le i\not= j\le \nu} A_i A_j) }{1-\E(\sum_{i=1}^\nu A_i^2)}.
	 $
By the dominated convergence theorem, we have that when $ \kappa\in (2, \infty]$, \begin{equation} \label{asymp-M>2} 
\E \frac{M_\infty^2}{a+M_\infty} \sim
\frac{1}{a}\, \E(M_\infty^2), \qquad a \to \infty.
\end{equation}

When $\kappa \in (1, 2]$, it is elementary to deduce from \eqref{cM}   that as $a \to \infty$, 
\begin{equation} \label{asymp-M}
\E \frac{M_\infty^2}{a+M_\infty} 
\, \asymp\,
  \begin{cases}
    a^{1- \kappa} ,
	& \hbox{\rm if $1 < \kappa < 2$} \, ,
    \\
   \frac{\log a}{a},    & \hbox{\rm if $ \kappa= 2$} .
    \end{cases}
\end{equation}


Recall from  \eqref{B2} that $\E \big( \frac{B_\varepsilon^2}{1+B_\varepsilon}\big)= \varepsilon\, \E\big( \frac1{1+B_\varepsilon}\big)$ which can be re-written as  \begin{equation}
	 \E \frac{\langle B_\varepsilon\rangle^2}{a+\langle B_\varepsilon\rangle} =  a \, \varepsilon\, \E \Big(  \frac1{1+B_\varepsilon}\Big) \sim a \, \varepsilon,   \qquad \varepsilon\to0, \label{aepsilon}
\end{equation}

\noindent where $a\equiv a(\varepsilon):= 1/\E(B_\varepsilon) \to \infty $ by \eqref{EB1/2}.  By Lemma \ref{l:phia} with $p=1$ there,    $$ \E \frac{\langle B_\varepsilon\rangle^2}{a+\langle B_\varepsilon\rangle} \le \E \frac{M_\infty^2}{a+M_\infty}.$$

\noindent Hence for $a = 1/\E(B_\varepsilon)$,  $$ \liminf_{\varepsilon\to0} \frac{1}{a\,  \varepsilon}\,  \E \frac{M_\infty^2}{a+M_\infty} \ge1,$$

\noindent which in view of \eqref{asymp-M} and \eqref{asymp-M>2} yield the Lemma. \hfill$\Box$

\medskip
 \medskip

We are ready to give  the asymptotic behaviors of $B_\varepsilon$ when $\kappa\in (2, \infty]$:

\begin{proposition} \label{p:kappa>2} Assume \eqref{hyp1} and \eqref{hyp2}.  If $\kappa \in (2, \infty]$, then under the probability $\P$, as $\varepsilon\to 0$, 
	$$\varepsilon^{-1/2} \, B_\varepsilon 
	\, \wcv \, 
	c_2\,  M_\infty,
	$$  
	with $c_2:=\big( \frac{1-\E(\sum_{i=1}^\nu A_i^2)}{\E(\sum_{1\le i\not= j\le \nu} A_i A_j)} \big)^{1/2}$ as in Proposition  \ref{p:B}. Moreover,  $$ 
	\lim_{\varepsilon\to0} \varepsilon^{-1/2} \, \E(B_\varepsilon) 
	= 
	c_2.$$
\end{proposition} 

{\noindent\bf Proof:} Based on  the boundedness in $L^2$ of $\varepsilon^{-1/2} B_\varepsilon$  (cf. \eqref{Bp1}),  it suffices to prove the convergence in law.  Let us first show the tightness of $ \frac{B_\varepsilon}{\sqrt \varepsilon} $ as $\varepsilon\to 0$.  By \eqref{EB1/2} and \eqref{kappa>2},  \begin{equation}\label{infsup}
	 c_2\le \liminf_{\varepsilon\to0} \frac{\E (B_\varepsilon)}{\varepsilon^{1/2}} 
	\le 
	\limsup_{\varepsilon\to0} \frac{\E (B_\varepsilon)}{\varepsilon^{1/2}} 
	\le 2. \end{equation}



In particular, under $\P$,  the family of the laws of  $(\frac{B_\varepsilon}{\sqrt \varepsilon}, \varepsilon \to 0)$ is tight. 
Take an arbitrary  subsequence $\varepsilon_n \to0$ such that $ \frac{B_{\varepsilon_n} }{\sqrt \varepsilon_n}\, \wcv \xi$,  with  some nonnegative random variable  $\xi$.  By \eqref{infsup},  $\xi$ is not degenerate; moreover, we deduce from \eqref{Blaw} that  $\xi$ must  satisfy  the cascade equation: $$
	\xi 
	\, \law\, 
	 \sum_{i=1}^\nu A_i \, \xi_i,
	 $$

\noindent where conditioned on $(A_i)$, $\xi_i$ are i.i.d. copies of $\xi$. By the uniqueness of the solution (see Liu \cite{liu00}),   $\xi= c\, M_\infty$ for some positive constant $c$.  We re-write \eqref{B2} as  $$ \E \Big( \frac{ (\frac{B_\varepsilon}{\sqrt \varepsilon})^2}{1+B_\varepsilon}\Big)=   \E\Big( \frac1{1+B_\varepsilon}\Big),$$

\noindent which by Fatou's lemma  along the  subsequence $\varepsilon_n\to0$, gives that $c^2 \, \E (M_\infty^2)\le 1$, i.e. $c\le  (\E (M_\infty^2))^{-1/2} =  c_2$.  This in view of the lower bound in \eqref{infsup} imply that $c= c_2$. Then we have proved that  any subsequence of $ \frac{B_\varepsilon}{\sqrt \varepsilon}$ converges to the same limit $c_2M_\infty$,  which  gives the Proposition. \hfill$\Box$

\section{Proof of Proposition \ref{p:B}}\label{s:4}

In view of Proposition \ref{p:kappa>2},   the proof of Proposition \ref{p:B} reduces  to show  the following two statements: As $\varepsilon\to0$, 
\begin{eqnarray}
	\E (B_\varepsilon(\varnothing))
	 &\asymp &  \begin{cases}
  \varepsilon^{1/\kappa},
	& \hbox{\rm if $1 < \kappa < 2$} \, ,
    \\
   \big(\frac{\varepsilon}{\log 1/\varepsilon}\big)^{1/2},    & \hbox{\rm if $ \kappa= 2$} \, . 
        \end{cases}
						\label{EB}
	\\  \nonumber
	\\
	\langle B_\varepsilon(\varnothing)\rangle \, 
	& \to & \, M_\infty, \qquad \hbox{for all $\kappa \in (1, \infty]$,  \, $\P$-a.s.} \label{B-proba}  	
\end{eqnarray}

 \noindent We remark that the $L^p$-convergence will follow from  the  $L^p$-boundedness of $\langle B_\varepsilon(\varnothing)\rangle $ given in \eqref{Bp1}.

Let us  start with a preliminary lemma which compares $\langle B_\varepsilon(\varnothing)\rangle $ with $M_\infty$. Write as before $B_\varepsilon\equiv B_\varepsilon(\varnothing)$ and denote by $\|\cdot \|_p:= \E(|\cdot |)^{1/p}$ the $L^p$-norm.

 \begin{lemma}\label{l:lp}  Assume \eqref{hyp1}, \eqref{hyp2} and let $p \in (1, \kappa) \cap (1, 2]$.    There exists some positive constant $c_{12}$ such that for any $\varepsilon \in (0, 1)$, 
 \begin{equation}\label{LpB-M}\| \langle B_\varepsilon\rangle - M_\infty\|_p   
 \le
  c_{12} \,  \times
 \begin{cases}
  (\E(B_\varepsilon))^{\frac\kappa{p}-1}, & \qquad \mbox{if $1<\kappa\le 2$}, 
  \\
 \varepsilon^{1/2}, &\qquad \mbox{if $\kappa \in (2, \infty]$}.
 \end{cases}
 \end{equation}
 \end{lemma}

 {\noindent\bf Proof:}   It follows from  \eqref{Bx}   that for any $x \in \T$, \begin{eqnarray}
	 \langle B_\varepsilon(x)\rangle 
	 &=&
	  \sum_{y: {\buildrel\leftarrow \over y}=x}  A(y)\,  \frac1{\E(B_\varepsilon)}\, \frac{\varepsilon+B_\varepsilon(y)}{1+B_\varepsilon(y)}  \nonumber
	  \\
	  &=& \sum_{y: {\buildrel\leftarrow \over y}=x}  A(y)\,    \langle B_\varepsilon(y)\rangle  +  \sum_{y: {\buildrel\leftarrow \over y}=x}  A(y)\,  \Delta(y), \label{Bx>}
\end{eqnarray}

\noindent with $$\Delta(y)
	:= \frac{\varepsilon}{\E(B_\varepsilon) (1+B_\varepsilon(y))}- 
	\frac1{\E(B_\varepsilon)}\, \frac{B_\varepsilon(y)^2}{1+B_\varepsilon(y)},
	$$

\noindent where  conditioned on $(A(y), {\buildrel\leftarrow \over y}=x, \nu_x)$,  $\Delta(y)$ are i.i.d. copies of \begin{equation}\label{Delta} \Delta:= \frac{\varepsilon}{\E(B_\varepsilon) (1+B_\varepsilon )}-  \frac1{\E(B_\varepsilon)}\, \frac{B_\varepsilon^2}{1+B_\varepsilon}. \end{equation}  We note from \eqref{B2} that $\E(\Delta)=0$.


  By iterating \eqref{Bx>}, we get that   for any $m\ge1$,  \begin{equation}\label{Bx>1}  \langle B_\varepsilon\rangle  \equiv \langle B_\varepsilon(\varnothing)\rangle  
  =
   \sum_{|x|=m}\, \prod_{\varnothing< y \le x} A(y)\,  \langle B_\varepsilon(x)\rangle + \Theta_m,\end{equation}

\noindent with   $$ \Theta_m:= \sum_{k=1}^m \sum_{|x|=k}\prod_{\varnothing< y \le x}  A(y) \Delta(x),$$

\noindent where conditioned on $(V(x), |x|\le m)$, $(B_\varepsilon(x), \Delta(x))$ are i.i.d. copies of $(B_\varepsilon , \Delta)$.    




Observe  that for any $m\ge 1$,  $M_\infty=   \sum_{|x|=m}\prod_{\varnothing< y \le x} A(y)\,  M_\infty^{(x)}$, where conditioned on $(A(x), |x|\le m)$, $M_\infty^{(x)}$ are i.i.d. copies of $M_\infty$.  Let  $$
	Y_m
	:= 
	\sum_{|x|=m} \prod_{\varnothing< y \le x} A(y)  \big( \langle B_\varepsilon(x)\rangle - M_\infty^{(x)}\big).$$

\noindent  Then we have \begin{equation}\label{M+Y+Theta} \langle B_\varepsilon \rangle = M_\infty + Y_m + \Theta_m.\end{equation}

To control $Y_m$, we use the following fact  (Petrov \cite{petrov}, pp. 82, (2.6.20)):  Let $k\ge1$ and $ 1\le p \le 2$. Let $\xi_1$, $\cdots$, $\xi_k$ be independent  random variables such that $\E(|\xi_i|^p)<\infty$ and $\E(\xi_i)=0$   
for all $ 1\le i\le k$. Then
\begin{equation}\label{petrov}
	 \E \left| \xi_1 + \cdots + \xi_k \right|^p 
	  \le
 	2 \sum_{i=1}^k \E(|\xi_i| ^p)  .
\end{equation}

\noindent Applying \eqref{petrov} to $Y_m$ yields  that  for any $p \in (1, \kappa) \cap (1, 2]$, $$\E \Big( |Y_m|^p\Big)
	  \le 
	2\, \E\Big(  \sum_{|x|=m}\prod_{\varnothing< y \le x} (A(y))^p \Big)\, \E \Big( |\langle B_\varepsilon \rangle- M_\infty|^p\Big) 
	=
	2\, a_p^m\, \E \Big( |\langle B_\varepsilon \rangle- M_\infty|^p\Big) ,$$
	
\noindent by using the constant $a_p\in (0, 1)$ given in \eqref{ap}.	Note that $\E(| \langle B_\varepsilon\rangle- M_\infty|^p)  \le  2^p \E(\langle B_\varepsilon\rangle^p) + 2^p \E (M_\infty^p) $ which is less than some positive constant $ c_{13}$  uniformly in  $\varepsilon \in (0, 1)$, by   \eqref{Bp1} and  \eqref{cM}.  It follows that  \begin{equation}
	 \E \Big( |Y_m|^p\Big)
	\le
	   2 c_{13}\, a_p^m, \qquad \forall m\ge1.  \label{Ymp}
\end{equation}



 By the triangular inequality,  $$  \|\Theta_m \|_p
\le
\sum_{k=1}^m \,   \E \Big( \big|\sum_{|x|=k}\prod_{\varnothing< y \le x}  A(y) \Delta(x)\big|^p\Big)^{1/p}.$$

\noindent For the above expectation term,  we apply \eqref{petrov}   and get  that \begin{equation}
	 \E \Big( \big|\sum_{|x|=k}\prod_{\varnothing< y \le x}  A(y) \Delta(x)\big|^p\Big)
	 \le 
	2\, \E\Big(  \sum_{|x|=k}\prod_{\varnothing< y \le x} (A(y))^p \Big)\, \E \Big( |\Delta|^p\Big) 
	=
	2\, a_p^k  \, \E \Big( |\Delta|^p\Big). \label{Thetam}\end{equation}

\noindent Consequently, for any $m\ge1$,  $$
\| \langle B_\varepsilon\rangle - M_\infty\|_p 
\le
\|Y_m\|_p + \|\Theta_m \|_p
\le
(2 c_{13})^{1/p} \, a_p^{m/p} + 2^{1/p} \sum_{k=1}^m  a_p^{k/p} \,   \| \Delta\|_p.$$

Recall  that  $0<a_p <1$ for $p \in (1, \kappa) \cap (1, 2]$. Letting $m\to \infty$ we get   that   \begin{equation}\label{B-M}
\| \langle B_\varepsilon\rangle - M_\infty\|_p 
\le
 \frac{2^{1/p}}{1- a_p^{1/p}}\,   \| \Delta\|_p, \qquad \forall \, \varepsilon \in (0, 1).
\end{equation}

 By \eqref{Delta}, $$ \| \Delta\|_p 
 \le
 \frac{\varepsilon}{\E(B_\varepsilon)} + \frac1{\E(B_\varepsilon)}\, \big \| \frac{B_\varepsilon^2}{1+B_\varepsilon} \big\|_p
 =
 a \varepsilon+ \big  \| \frac{ \langle B_\varepsilon\rangle ^2}{a + \langle B_\varepsilon\rangle} \big\|_p,
 $$

\noindent with $a\equiv a(\varepsilon):= \frac1{\E(B_\varepsilon)}$ as in \eqref{aepsilon}.   By   Lemma \ref{l:phia}, $$\E \Big(\frac{ \langle B_\varepsilon\rangle ^2}{a + \langle B_\varepsilon\rangle}\Big)^p
\le
\E \Big(\frac{ M_\infty^2}{a + M_\infty}\Big)^p, \qquad \forall\, \varepsilon \in (0, 1).$$

\noindent  When $\kappa \in (1, 2]$,  we deduce from \eqref{cM} that there exists some positive constant $c_{14}$ such that for all  $a>1$,  $$\E \Big(\frac{ M_\infty^2}{a + M_\infty}\Big)^p
\le c_{14}  \,  a^{-(\kappa-p)} .$$

\noindent When $\kappa \in (2, \infty]$, $\E(M_\infty^2)< \infty$ and we have  $\E  (\frac{ M_\infty^2}{a + M_\infty} )^p \le a^{-p} \E(M_\infty^p).$ It follows that for some positive constant $c_{15}$, $$
\big  \| \frac{ \langle B_\varepsilon\rangle ^2}{a + \langle B_\varepsilon\rangle} \big\|_p
 \le
 c_{15}  \times
 \begin{cases}
  (\E(B_\varepsilon))^{\frac\kappa{p}-1}, & \qquad \mbox{when $1<\kappa\le 2$}, 
  \\
  \E(B_\varepsilon), &\qquad \mbox{when $\kappa \in (2, \infty]$}.
 \end{cases}
$$

 Consider at first  the case $1 < \kappa < 2$.     By  \eqref{kappa<2}, there exists some positive constant $c_{16}$ such that   $\E(B_\varepsilon) \ge c_{16}\, \varepsilon^{1/\kappa}$, hence $$\| \Delta\|_p
  \le 
\frac1{c_{16}} \, \varepsilon^{1-1/\kappa} + c_{15}\,  (\E(B_\varepsilon))^{\frac\kappa{p}-1}
 \le
  c_{17} (\E(B_\varepsilon))^{\frac\kappa{p}-1}, \qquad \forall \, \varepsilon \in (0, 1),
   $$
 
 \noindent for some positive constant $c_{17}$. 
 
 For the case $\kappa=2$,  we use \eqref{kappa=2} and  the same argument as above to get that     $\| \Delta\|_p  \le c_{17} (\E(B_\varepsilon))^{\frac\kappa{p}-1}$ by eventually enlarging the constant $c_{17}$. 
 
 When $\kappa >2$, we have already proven in Proposition \ref{p:kappa>2} that $\E(B_\varepsilon) \asymp \varepsilon^{1/2}$, which yields that $\| \Delta\|_p \le a \varepsilon + c_{15} \E(B_\varepsilon) \le c_{17} \varepsilon^{1/2}$. 
 
 Finally  Lemma \ref{l:lp} follows from the above estimates on $\| \Delta \|_p$ and  \eqref{B-M}.   \hfill$\Box$

\medskip
We are now ready to give the proofs of  \eqref{EB} and \eqref{B-proba}:

{\noindent\bf Proof of  \eqref{EB}.} In this proof we treat the cases   $\kappa \in (1, 2]$.  
 By Lemma \ref{l:low},   we only need to prove the upper bounds of $\E(B_\varepsilon)$ in \eqref{EB}.

For any $r>1$, we have $$ \P  ( \langle B_\varepsilon\rangle > r  ) \ge 
\P(M_\infty >  2 r) - \P ( B_\varepsilon - M_\infty  \le -   r ).$$

\noindent Let $p\in (1, \kappa)$.   By Markov's inequality, $$\P ( B_\varepsilon - M_\infty  \le -   r ) \le 
r^{-p} \, \|B_\varepsilon - M_\infty \|_p^p 
\le  c_{12}^p \, r^{-p} \,  \E(B_\varepsilon)^{\kappa -p}   ,$$

\noindent where the last inequality follows from  Lemma \ref{l:lp}. 
By  \eqref{cM}, there exists some positive constant $c_*$ such that $$ 
	\P  ( M_\infty>  2 r ) 
	\ge  
	c_*\,  r^{-\kappa}, \qquad \forall\, r\ge1.
	$$

It follows that  for any $r \ge 1$, $$\P  ( \langle B_\varepsilon\rangle > r  ) \ge
c_* r^{-\kappa} -  c_{12}^p  \, r^{-p} \,  \E(B_\varepsilon)^{\kappa -p}  ,$$

\noindent which is larger than 	$c_* r^{-\kappa} /2$ if   furthermore $r $ is such that $c_* r^{-\kappa}  \ge 2 c_{12}^p \, r^{-p} \,  \E(B_\varepsilon)^{\kappa -p}  $, namely if $ r \le c_{18}   \E(B_\varepsilon)^{-1} , $  with $c_{18}:= \big(\frac{c_*}{2 c_{12}^p}\big)^{1/(\kappa-p)}$. In other words,  we have proved that  when $\kappa \in (1, 2]$, for any $1\le r \le c_{18}  \E(B_\varepsilon)^{-1}$,  \begin{equation} \label{probaB>r} \P \big( \langle B_\varepsilon\rangle > r \big) \ge
\frac{c_*}{2} r^{-\kappa} . \end{equation}
	
	 Recall from \eqref{aepsilon} that   \begin{equation}  \E \frac{\langle B_\varepsilon\rangle^2}{\frac1{\E(B_\varepsilon)}+ \langle B_\varepsilon\rangle}
 	 	 \sim 
	 \frac\varepsilon{\E(B_\varepsilon)}, \qquad \varepsilon\to0. 	\label{up4}
	\end{equation}

	Let $r:=   c_{18}  \E(B_\varepsilon)^{-1}$ which is larger than $1$ for all small $\varepsilon$. We have  $$ \E \frac{\langle B_\varepsilon\rangle^2}{\frac1{\E(B_\varepsilon)}+ \langle B_\varepsilon\rangle} 
	\ge 
	\frac{r^2}{\frac1{\E(B_\varepsilon)}+ r} \P \Big( \langle B_\varepsilon\rangle > r \Big) 
	\ge 
	c_{19}\,  \E(B_\varepsilon)^{\kappa -1},
$$

\noindent  for some positive constant $c_{19}$ independent of $\varepsilon$.  This in view of \eqref{up4} imply   that when $\kappa \in (1, 2],$  there exists some positive constant $c_{20}$ such that  \begin{equation}\label{uppEB}  \E (B_\varepsilon) \le c_{20}\, \varepsilon^{1/\kappa} , \qquad \forall \, \varepsilon \in (0, 1), \end{equation}

\noindent which gives the desired upper bound in \eqref{EB} when $\kappa\in (1, 2)$.  

To deal with  the case $\kappa=2$, we write as before  $a: =  1/\E(B_\varepsilon)$. By integration by parts,  we deduce from \eqref{probaB>r} that for all sufficiently small $\varepsilon$ such that $a   >\max(1, 1/c_{18})$, \begin{eqnarray*} \E \frac{\langle B_\varepsilon\rangle^2}{\frac1{\E(B_\varepsilon)}+ \langle B_\varepsilon\rangle}
& = &
 \int_0^\infty \frac{x(2a +x)}{(a+x)^2} \P(\langle B_\varepsilon\rangle> x) d x 
 \\
 &\ge&
 \int_1^{c_{18} \, a}  \frac{x(2a +x)}{(a+x)^2} \P(\langle B_\varepsilon\rangle> x) d x 
\\
& \ge &   \frac{c_*}{2} \int_1^{c_{18} \, a} \frac{x(2a +x)}{(a+x)^2} x^{-2} d x
\\
 &\ge& 
c_{21} \, \frac{\log a}{a}, \end{eqnarray*} 

\noindent  for some positive constant $c_{21}$ independent of $\varepsilon$.  It follows from \eqref{up4} that  there is some positive constant $c_{22}$ such that  $a \ge c_{22} \sqrt{\frac{\log 1/\varepsilon}{\varepsilon}}$ for all $\varepsilon \in (0, 1)$. This gives the desired upper bound  in \eqref{EB} for the case $\kappa=2$ and completes the proof of \eqref{EB}.  \hfill$\Box$

\medskip 

Now  to prove Proposition \ref{p:B}, it remains to show the almost sure convergence in  \eqref{B-proba}.
 
{\noindent\bf Proof of  \eqref{B-proba}.}   By \eqref{LpB-M},  for any  $\kappa\in (1, \infty]$ and $p\in (1, \kappa) \cap (1, 2]$, there exists some positive constant $\varrho$ such that   
 \begin{equation}\label{LpB-M2}\| \langle B_\varepsilon \rangle - M_\infty\|_p   \le  c_{12} \,   \varepsilon^{\varrho}, \qquad 0< \varepsilon <1.
 \end{equation}

  Let $\varepsilon_n:= n^{-2/\varrho}$.  It follows from \eqref{LpB-M2} that $$ \sum_{n=1}^\infty \| \langle B_{\varepsilon_n} \rangle- M_\infty\|_p < \infty,$$

\noindent which yields that $\P$-a.s., $ \langle B_{\varepsilon_n} \rangle \to M_\infty$ as $n \to \infty$.   Observe that $\varepsilon\to B_\varepsilon$ is non-increasing, hence for any $\varepsilon_n \le \varepsilon < \varepsilon_{n-1}$, $ \langle B_{\varepsilon_{n-1}} \rangle \, \frac{\E(B_{\varepsilon_{n}})}{\E(B_{\varepsilon_{n-1}})}\le  \langle B_\varepsilon  \rangle \le \langle B_{\varepsilon_n} \rangle \, \frac{\E(B_{\varepsilon_{n-1}})}{\E(B_{\varepsilon_n})}$, then \eqref{B-proba} follows immediately if we can show  that  $\frac{\E(B_{\varepsilon_{n-1}})}{\E(B_{\varepsilon_n})} \to 1$ as $n \to \infty$. 

To this end, 
define  $b(\varepsilon):= \E(B_\varepsilon)$ for $0< \varepsilon<1$ and $b(0):=0$. By \eqref{Bepsilon}, we have $$ b(\varepsilon)=
\E \Big[ 1- E_{\varnothing, \omega}\big((1-\varepsilon)^{\frac12 (1+T_{{\buildrel\leftarrow \over \varnothing}})}\big)\Big].$$

\noindent Then $$b'(\varepsilon)= \E\otimes E_{\varnothing, \omega} \Big[ \frac12 (1+T_{{\buildrel\leftarrow \over \varnothing}}) (1-\varepsilon)^{\frac12 (T_{{\buildrel\leftarrow \over \varnothing}}-1)} \Big].$$

\noindent Since $T_{{\buildrel\leftarrow \over \varnothing}}\ge 1$,  $P_{\varnothing, \omega}$-a.s., we see that $b'$ is a positive decreasing function  on $(0, 1)$.   By \eqref{EB} and Proposition \ref{p:kappa>2}, there exists some constant  $c_{23}>1$ such that $$ \frac1{c_{23}}  r(\varepsilon) \le f(\varepsilon) \le c_{23} r(\varepsilon),$$

\noindent where $r(\varepsilon):= \varepsilon^{1/\kappa}$ when $\kappa \in (1, 2)$,   $r(\varepsilon):= \sqrt{\varepsilon/\log(\ee/\varepsilon)}$ when $\kappa=2$ and $r(\varepsilon):=\varepsilon^{1/2}$ when $\kappa \in (2, \infty]$.  By the concavity of $b$,  $$b'(\varepsilon) \le \frac{b(\varepsilon)}{\varepsilon} \le c_{23} \frac{r(\varepsilon)}{\varepsilon}.$$


\noindent It follows that  $$ \frac{b(\varepsilon_{n-1})- b(\varepsilon_n)}{b(\varepsilon_n)}  \le (\varepsilon_{n-1}- \varepsilon_n) \frac{b'(\varepsilon_n)}{b(\varepsilon_n)} \le
  (\varepsilon_{n-1}- \varepsilon_n)  \frac{c_{23}^2}{\varepsilon_n}  \to 0,  \qquad n \to\infty, $$ 

\noindent  by the choice of $\varepsilon_n$. Hence $\frac{\E(B_{\varepsilon_{n-1}})}{\E(B_{\varepsilon_n})} \to 1$ as $n \to \infty$,  we get    \eqref{B-proba} and complete   the proof of Proposition \ref{p:B}. \hfill$\Box$

\bigskip
{\noindent\bf Acknowledgements.} I am grateful to Pierre Rousselin who has kindly pointed out a mistake in the original version.

\end{document}